# IDENTIFICATION OF TARGET SYSTEM OPERATIONS. DEVELOPMENT OF GLOBAL EFFICIENCY CRITERION OF TARGET OPERATIONS

**I. Lutsenko**
PhD, Professor
Department of Electronic Devices
Kremenchuk Mykhailo Ostrohradshyi
National University
Pervomaiskaya str., 20,
Kremenchuk, Ukraine, 39600
E-mail: delo-do@i.ua

*Розроблено глобальний показник ефективності цільових операцій (критерій оптимального управління). Встановлено, що даний показник дозволяє дати системну оцінку досліджуваним операціям. Більш високе значення ефективності конкретної операції, по відношенню до порівнюваної, говорить про те, що дана операція забезпечить накопичення більш високої доданої вартості (прибутку) з урахуванням витрат і часу операції*

*Ключові слова: ефективність операції, формула ефективності, ефективність використання ресурсів, критерій ефективності, критерій оптимізації*

*Разработан глобальный показатель эффективности целевых операций (критерий оптимального управления). Установлено, что данный показатель позволяет дать системную оценку исследуемым операциям. Более высокое значение эффективности исследуемой операции, по отношению к сравниваемой операции, говорит о том, что данная операция обеспечит накопление более высокой добавленной стоимости (прибыли) с учетом затрат и времени операции*

*Ключевые слова: эффект, эффективность, эффективность использования ресурсов, критерий эффективности, формула эффективности, критерий оптимизации*

## 1. Introduction

The ability to consciously affect the financial result of the economic system is the main goal of any qualified and properly motivated management. One of the major tools of such influence is the ability of unique identification of the operations studied in order to select the most effective one.

Currently, there is quite a strange situation in the knowledge "market". A great number of publications consider efficiency as the "old good friend". At the same time, without the validation of the developed criterion, a reader is invited to use one or the other indicator as efficiency indicator. It resembles the actions of healer who offers the patient a miracle drug in the hope of a placebo effect.

The result of this approach can be seen everywhere. Automation of technological processes has reached an impressive level, unlike the automation of management processes.

*What is the reason?*

The reason is that it is very risky to commit financial decision–making to automation. Because the lack of an adequate optimization criterion (real efficiency indicator) can automatically quickly cause significant financial damage to the owner.

Therefore, for example, the huge army of automatic systems operates in management modes with strongly reduced functionality. The most important functions that should provide the maximum economic effect of a functioning system do not work there. Hence, the development of a single cybernetic efficiency indicator will allow to raise management process automation issues to a new level.

## 2. Analysis of published data and problem statement

In operations research, some of the evaluation basics of target operations were laid in the first twenty years after the creation of this discipline. Since then, there was no significant progress in this regard [1]. Herewith, it is natural that the focus is more shifted towards the use of global civilization opportunities, such as GPS–navigation systems [2], remote data collection using the Internet technologies, mathematical modeling methods [3]. Fundamental changes are not observed in inventory management problems [4], or optimal enterprise management issues [5].

In the overwhelming majority of works, which are directly connected with the estimation theory, the issue of the optimization criterion adequacy is not even raised. Although the problem has long been known for decades.

Attempts to solve the efficiency formula derivation problem throughout the development history of the management theory, economic theory, operations research, etc. have shown that traditional approaches based on the fundamental theories that underlie yesterday's ideas of cybernetics, have failed: "Engineers, researchers, economists and designers are continuously suggesting "universal, accurate and clear" objective functions. In 1967, one of the authors (L. Barsky) managed to collect over a hundred optimization criteria of separation processes [6]. Their classification has revealed that there is no universal criterion, and the selection of process optimization or efficiency criterion is not an easy task" [7].

These failures are caused by an attempt to use the traditional approach to the creation of a model of the target operation, the foundation of which was laid not in cybernetics, but in the economy. Herewith, there is no theoretical justification of the basic model of the economic operation at all.

Creation of a fundamentally new concept, model of the deployed target operation has allowed to lay the foundation for moving forward in this direction. Determination of the time of actual completion of the target operation [8] and its resource intensity [9] made it possible to identify the operations in terms of their efficiency for a number of special cases.

However, for the majority of optimal management, operations research and economic analysis problems, it is necessary to assess target operations with arbitrarily chosen or set management parameters.

This is a prerequisite for using the efficiency criterion at all process optimization stages.

### 3. The goal and objectives of the study

The goal is to develop a unified and reliable formula of the resource efficiency.

For this purpose, the following tasks were solved:
– development of test operations to validate the efficiency indicator;
– determination of the relationship between the theory of efficiency and economic profitability;
– determination of the indicator of the potential effect of the target operation;
– determination of the efficiency indicator;
– verification of the efficiency indicator during the study of test operations.

### 4. Test operations to validate the efficiency indicator

Before proceeding to the efficiency indicator development, it is necessary to answer the question of how to evaluate the result?

One way is to estimate the special reference operations with the same or very similar efficiency. Four groups of such simple reduced operations are shown in Fig. 1.

*What is characteristic of these operations?*

These operations are characterized by the fact that the subtotal of these operations has the same value to a certain point in time. For example, all operations with the subtotal at the time $t=12$ have the cost estimate of released – invested resources at the level of 5,315 monetary units.

The operations in each group have their duration in time. For the first group of operations this time is $T_{op.1}=2$ time intervals (t. int.), for the second group – $T_{op.2}=4$ (t. int.), the third – $T_{op.3}=6$ (t. int.), and the fourth – $T_{op.4}=8$ (t. int.).

The initial cost estimate of financial resources invested in the first operation of each group is the same – $RE_{1.1} = RE_{2.1} = RE_{3.1} = RE_{4.1} = 3$ monetary units.

The cost estimate of financial resources invested in each subsequent operation of the group is equal to an amount of released funds from the previous operation.

| RVIC | | | | Time |
|---|---|---|---|---|
| 1,1 | 1,21 | 1,331 | 1,4641 | |
| Reference operations | | | | |
| Operations of the 1 group | Operations of the 2 group | Operations of the 3 group | Operations of the 4 group | |
| 3 | 3 | 3 | 3 | 0 |
| | | | | 1 |
| 3,3 | | | | 2 |
| | | | | 3 |
| 3,63 | 3,63 | | | 4 |
| | | | | 5 |
| 3,993 | | 3,993 | | 6 |
| | | | | 7 |
| 4,392 | 4,392 | | 4,392 | 8 |
| | | | | 9 |
| 4,832 | | | | 10 |
| | | | | 11 |
| 5,315 | 5,315 | 5,315 | | 12 |
| | | | | 13 |
| 5,846 | | | | 14 |
| | | | | 15 |
| 6,431 | 6,431 | | 6,431 | 16 |

Fig. 1. Four groups of simple target reference operations

Thus, the released financial resources of the first operation of the second group in an amount of $RE_{2.2} = PE_{2.1} = 3{,}63$ monetary units are the initial cost estimate of the funds invested in the second operation of this group, etc. That is, the value added is capitalized.

The amount of the released funds of the operation with respect to the invested funds in each group is proportional to the resource value increase coefficient (RVIC) of the target operation, $k_1$, $k_2$, $k_3$ and $k_4$. For the second group of operations, the value of this coefficient is $k_2=1.21$.

Accordingly, the cost estimate of invested and released funds, for example, of the operations of the second group is defined by the expression $PE_2 = k_2 \cdot RE_2$.

RVICs are selected so that the operations that end up at the same time have the same amount of funds released. Since the initial cost estimate of financial resources invested in the first operation of each group is the same, so it is at the time of their simultaneous release, from the point of view of resource efficiency, all these operations, to a first approximation, can be considered equivalent.

Thus, the integral indicator, which links the three basic indicators RE, PE, $T_{op}$, and thus identifies each operation of the Fig. 1 as equivalent, is the efficiency indicator.

### 5. The efficiency of the target operation

Previously obtained "resource intensity" indicator [5] shows the losses of management of the target operation. The value added of this operation (let's define it as a parent operation) can be converted into resources and used to conduct an independent operation, which we define as a daughter operation.

We define the ratio of the resource intensity of the daughter operation of the parent operation as efficiency of products transformation process.

Investments in future operation are determined by the quantity of value added of the previous operation. If the cost estimate of invested funds of the previous operation is determined by the RE, and the cost estimate of released funds – by PE, then the volume of funds PE−RE will be invested in the next operation. Then the value added of future operation can be defined by the expression k (PE−RE).

By substituting the values obtained in the expression for determining the resource intensity of the simple reduced operation, we obtain the value of the resource intensity of the daughter operation.

$$R_d = \frac{k(PE-RE)(PE-RE)T_1^2}{2[k(PE-RE)-(PE-RE)]} = \frac{kT_d^2}{2(k-1)}(PE-RE),$$

$$k > 1,$$

where $R_d$ – resource intensity of the daughter operation; $k$ – resource value increase coefficient of the previous operation; $T_d$ – time interval of the daughter operations;

$$RE = \int_{t_0}^{t_a}\left(\int_{v_0}^{v}|re(v)|dv\right)dt, \ v \in [0, \ t_a];$$

$$PE = \int_{t_0}^{t_a}\left(\int_{v_0}^{v}pe(v)dv\right)dt, \ v \in [0, \ t_a];$$

$t_p$ – time of the actual completion of the parent target operation.

Let us define the efficiency of the target operation ( Eff ) as the ratio of resource intensity of the daughter operation to the resource intensity of the parent operation

$$Eff = \frac{\frac{kT_d^2}{2(k-1)}(PE-RE)}{\frac{PE \cdot RE \cdot (T_p)^2}{2 \cdot (PE-RE)}} = \frac{\frac{kT_d^2}{(k-1)}(PE-RE)^2}{PE \cdot RE \cdot (T_p)^2}, \ k > 1.$$

Let us assume that the coefficient of the value added does not change from one operation to another, that is $PE = k \cdot RE$. Then, for the target operations with a unit duration in time we obtain

$$Eff = \frac{k \cdot RE - RE}{RE} = \frac{PE - RE}{RE}, \ PE > RE.$$

This expression is nothing but the economic profitability of the operation. Consequently, the profitability is a particular case of the efficiency of a pair of target operations with unit time intervals.

## 6. The potential effect of the target operation

Let us return to the model of the target operation [4]. At the time $t_a$, the thread $ide(\tau)$ stops compensating the closed thread of tight resources $ibe(\tau)$ and begins generating the effect of the target operation (Fig. 2). The magnitude of this effect is determined by the area, which is limited to the function $vde(\tau)$.

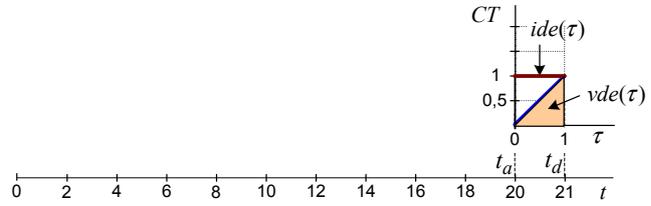

Fig. 2. Illustration of the function $vde(\tau)$

Since the magnitude of the future effect is beyond the management of the operation studied, only the potential effect within the unit time interval can be assessed. The magnitude of this effect is determined by the value of the function $wde(\tau)$ at a time $t_a + 1 = t_d$ (Fig. 3).

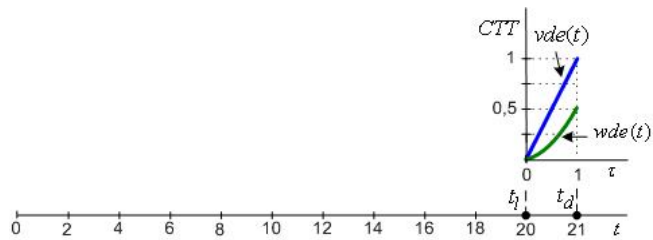

Fig. 3. Illustration of the function $wde(\tau)$

On the example of the registration operation, let us consider the steps to be taken for determining the absolute potential effect of the target operation.

1. We determine the TACO from the expression

$$t_a = \frac{PE \cdot t_p - RE \cdot t_r}{PE - RE}.$$

2. We determine the position of the right boundary point of the function domain

$$t_d = t_a + \Delta t = t_a + 1.$$

3. We introduce auxiliary variables $v$, $\tau$ and define them

$$v \in [t_0; t_d]; \ \tau \in [t_a; t_d].$$

4. We build a deployed model of the operation

$$ice(v) = \int_0^v pe(s)ds + \int_0^v re(s)ds.$$

5. We determine the function $vde(\tau)$ on the interval $\tau \in [t_a; t_d]$ (Fig. 2)

$$vde(\tau) = \int_{t_a}^{\tau}\left[\int_{t_0}^{v}pe(s)ds + \int_{t_0}^{v}re(s)ds\right]dv.$$

6. We determine the function $wde(\tau)$ on the interval $[t_l; t_d]$ (Fig. 3)

$$wde(t) = \int_{t_a}^{t} \left( \int_{t_a}^{\tau} \left[ \int_{t_0}^{v} pe(s)ds + \int_{t_0}^{v} re(s)ds \right] dv \right) d\tau.$$

The potential effect of the IS operation is numerically equal to the value of the function $wde(\tau)$ at the time $t_d$. Therefore, the expression for determining the absolute potential effect of the operation in case of a vector representation of the original model, will have the form

$$A = \int_{t_a}^{t_d} \left( \int_{t_a}^{\tau} \left[ \int_{t_0}^{v} pe(s)ds + \int_{t_0}^{v} re(s)ds \right] dv \right) d\tau,$$

where $A$ – the potential effect of the studied target operation.

Based on the geometric interpretation of the definition of the absolute potential effect (APE) indicator, its numerical value can be determined from the expression

$$A = \left[ (PE - RE)(t_d - t_l)^2 \right] / 2 = \left[ (PE - RE) T_l^2 \right] / 2.$$

Thus, the absolute potential effect of the target operation is determined by the half of the product of the value added (cost) by the square of the estimated time interval.

### 7. The expression for determining the resource efficiency of the target operation

By determining the magnitude of the potential effect of the studied operation, we can obtain an expression for determining the efficiency, in general terms, as the ratio of the potential effect of the studied operation to its resource intensity

$$E = \frac{A}{R},$$

where $E$ – the efficiency of the studied target operation.

When determining the target operations with distributed parameters, expression for determining the efficiency of the studied operation will be of the form

$$E = \frac{\int_{t_a}^{t_d} \left( \int_{t_a}^{\tau} \left[ \int_{t_0}^{v} pe(s)ds + \int_{t_0}^{v} re(s)ds \right] dv \right) d\tau}{\int_{t_0}^{t_a} \left[ \int_{v_0}^{v} \left| \int_{v_0}^{v} re(v)dv \right| dv - \int_{v_0}^{v} \left( \int_{v_0}^{v} pe(v)dv \right) dv \right] dv},$$

$v \in [0, t_a]$, $k > 1$.

For evaluating models of simple target operations, it is enough to use the simplified analytical expression

$$E = \frac{(PE - RE)^2 T_l^2}{PE \cdot RE \cdot T_{op}^2}, \quad k > 1.$$

Let us use this expression to solve the problem associated with the estimation of reference operations (Table 1). But first, we let us perform a little comparative study of the obtained cybernetic criterion of "efficiency" with the global economic indicator such as "profitability".

Each target operation of a set has the cost estimate of input products of the operation (RE – costs), the cost estimate of output products of the operation (PE) and the operation time ($T_{op}$). For each operation, efficiency (E) that is considered on the diagrams in conjunction with profitability (PROF) is calculated.

Let us consider the way the profitability and efficiency respond to a change in parameters of target operations in those cases when these changes are evident and predictable.

The first set of operations (Table 1) is characterized by the fact that from one operation to another, the cost estimate of input products of the operation RE increases and the cost estimate of output products of the operation PE and the operation time $T_{op}$ do not increase.

Table 1

A set of registration models of target operations for investigating the relationship of efficiency and profitability when changing the cost estimate of input products

| N | RE | PE | $T_{op}$ | PROF | E | k |
|---|----|----|----|------|----|----|
| 1 | 2   | 3 | 1 | 0,50 | 0,50 | 1,5 |
| 2 | 2,1 | 3 | 1 | 0,43 | 0,43 | 1,429 |
| 3 | 2,2 | 3 | 1 | 0,36 | 0,36 | 1,364 |
| 4 | 2,3 | 3 | 1 | 0,30 | 0,30 | 1,304 |
| 5 | 2,4 | 3 | 1 | 0,25 | 0,25 | 1,25 |
| 6 | 2,5 | 3 | 1 | 0,20 | 0,20 | 1,2 |
| 7 | 2,6 | 3 | 1 | 0,15 | 0,15 | 1,154 |

In this case, with an increase in the cost estimate of input products of the operation (costs RE), at constant PE and $T_{op}$, efficiency of the operation should reduce. Calculation of efficiency for the first set of operations (Table 1) confirms this assumption (Fig. 4). It can be seen that the efficiency has higher dynamics of the decrease with respect to profitability.

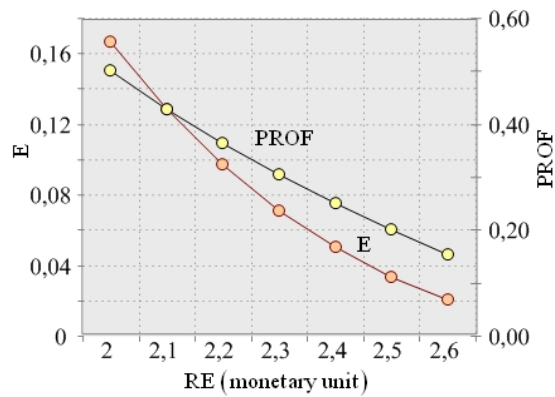

Fig. 4. Changes in the efficiency E and profitability PROF from the costs RE of the target operation

The second set of operations (Table 2) is characterized by the fact that from one operation to another, the cost estimate of output products of the operation PE increases and the cost estimate of input products of the operation RE and the operation time $T_{op}$ do not change. The higher the cost estimate of the output product of the target operation, the higher its effectiveness should be.

Table 2

A set of registration models of target operations for investigating the relationship of efficiency and profitability when changing the cost estimate of output products

| N | RE | PE | $T_{op}$ | PROF | E | k |
|---|----|----|----|------|---|---|
| 1 | 2 | 2,5 | 1 | 0,25 | 0,25 | 1,25 |
| 2 | 2 | 2,6 | 1 | 0,30 | 0,30 | 1,3 |
| 3 | 2 | 2,7 | 1 | 0,35 | 0,35 | 1,35 |
| 4 | 2 | 2,8 | 1 | 0,40 | 0,40 | 1,4 |
| 5 | 2 | 2,9 | 1 | 0,45 | 0,45 | 1,45 |
| 6 | 2 | 3 | 1 | 0,50 | 0,50 | 1,5 |
| 7 | 2 | 3,1 | 1 | 0,55 | 0,55 | 1,55 |

Calculation of efficiency for the second set of operations (Table 2) confirms this hypothesis as well (Fig. 5). The growth rate of efficiency is slightly higher than that of profitability.

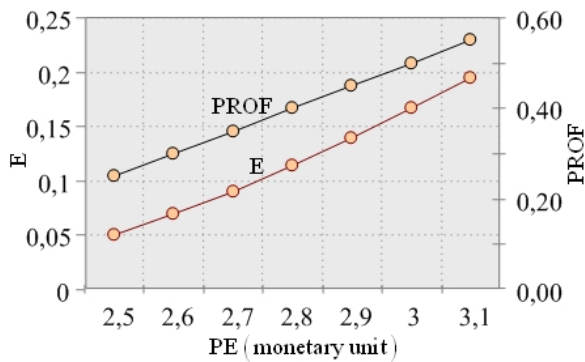

Fig. 5. Change in the efficiency E and profitability PROF from the expert (cost) estimate of output products PE of the target operation

In the third set of operations, values RE and PE do not change, and the operation time $T_{op}$ changes (Table 3).

Table 3

A set of registration models of target operations for investigating the relationship of efficiency and profitability when changing the target operation time

| N | RE | PE | $T_{op}$ | PROF | E | k |
|---|----|----|------|------|---|---|
| 1 | 2 | 3 | 1 | 0,5 | 0,5 | 1,5 |
| 2 | 2 | 3 | 1,2 | 0,5 | 0,35 | 1,5 |
| 3 | 2 | 3 | 1,3 | 0,5 | 0,3 | 1,5 |
| 4 | 2 | 3 | 1,4 | 0,5 | 0,26 | 1,5 |
| 5 | 2 | 3 | 1,5 | 0,5 | 0,22 | 1,5 |
| 6 | 2 | 3 | 1,6 | 0,5 | 0,2 | 1,5 |
| 7 | 2 | 3 | 1,7 | 0,5 | 0,17 | 1,5 |

Obviously, the longer the operation time (at fixed RE and PE), the longer the input products RE of the operation are bound by technological processes, and the lower the efficiency should be.

As can be seen, the "efficiency" indicator shows a tendency to decline (Fig. 6).

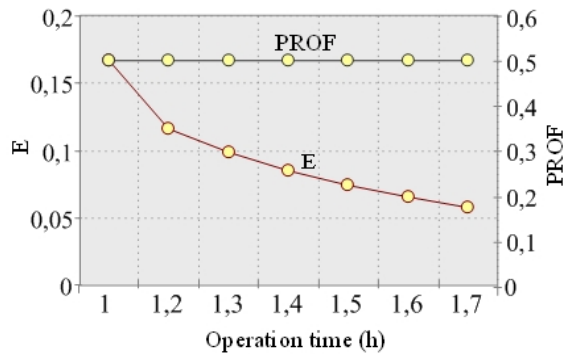

Fig. 6. Change in the efficiency E from the target operation time

The profitability "does not solve" such problems in principle.

Table 4 shows the results of evaluating the efficiency for the groups of models of operations shown in Fig. 1.

Table 4

The efficiency of groups of target operations (Fig. 1)

| Operations | E |
|---|---|
| 1 | 0,002273 |
| 2 | 0,002278 |
| 3 | 0,002287 |
| 4 | 0,002298 |

A natural question that may arise in the study of data – what is the reason for the small but steady growth of efficiency while an increase in the operations group number? (Table 4).

The reason for this growth is that the operations in the Table actually are not equivalent in terms of equality of their efficiencies. Their efficiencies are close but not equal.

For example, the operations in the first group have a lower efficiency compared to the operations in the second group. In order to deal with the cause, it is possible, for example, to represent the first two operations of the first group as three operations (Fig. 7). It is seen that the equivalent circuit of the lower group of operations requires additional investments in an amount of 0.3 monetary units. Therefore, the efficiency of the first operation of the second group is really higher.

All of this suggests that we can not "rigidly" set the parameters of all three variables when constructing reference operations.

Now when we have a "candidate" for the efficiency formula, it is possible to determine the parameters of the reference target operation with respect to efficiency if the basic reference operation is known.

Let the basic reference operation be the first operation of the first group. Let us determine its efficiency:

$$E_{1.1} = \frac{(PE-RE)^2 T_1^2}{PE \cdot RE \cdot T_{op}^2} = \frac{(3,3-3)^2 1^2}{3,3 \cdot 3 \cdot 2^2} = 0,002272727.$$

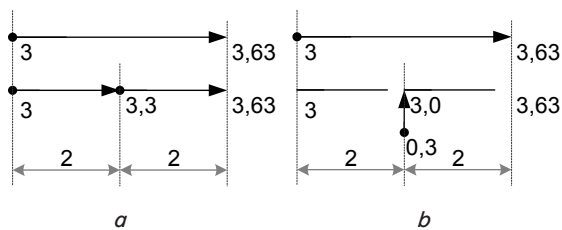

Fig. 7. Circuit of the group of operations: *a* — original representation of the first two operations of the first group with respect to the first operation of the second group; *b* — representation of the operations of the first group in the form of the equivalent circuit of the three operations

For the first operation of the second group to be equivalent with respect to the efficiency of the basic operation, it is necessary to make one of the parameters of this operation free. Let it be PE. By calculating PE with the given parameters $RE = 3$, $T_{op} = 4$ and $E = 0{,}002272727$, we obtain $PE = 3{,}6292174$.

Target operations in the group are equivalent. For example, let us compare the first and second operations of the fourth group

$$E_{2.1} = \frac{(PE - RE)^2 T_1^2}{PE \cdot RE \cdot T_{op}^2} = \frac{(3{,}63 - 3)^2 1^2}{3{,}63 \cdot 3 \cdot 8^2} = 0{,}000569473,$$

$$E_{2.2} = \frac{(PE - RE)^2 T_1^2}{PE \cdot RE \cdot T_{op}^2} = \frac{(4{,}3923 - 3{,}63)^2 1^2}{4{,}3923 \cdot 3{,}63 \cdot 8^2} = 0{,}000569473.$$

It should be noted that the efficiency formula is not a magic elixir. For its practical use in manufacturing enterprises, it is necessary to change approaches to the creation of technological cycles, and managed systems.

But the possibility of using a single optimal management criterion at all, without exception, stages of the technological and trading processes at least makes this work real and feasible in the near future.

## 8. Conclusions

The method of creating a set of operations that can be used as a test for preliminary evaluation of the adequacy of the indicators, claiming to be the efficiency criterion was developed. The method lies in forming strings of successive operations with the same time and the same RVIC of each operation of a specific chain. Herewith, the cost estimate of input products of each subsequent operation is equal to the cost estimate of output products of the previous operation. Moreover, the operation time in each chain is set by the multiple of time of the shortest operation of the entire set of operations, and RVIC such that at the time of simultaneous completion of the operations in different chains, cost estimates of their output products are the same.

It was found that the expression that defines the efficiency of the target operation is equivalent to the expression of economic profitability for a set of target operations with unit time interval, since, in this case, these expressions are identical.

The definition of the cybernetic term of "potential effect" of the operation was given and the possibility for its quantification was revealed. The absolute potential effect is determined by the integral evaluation of the integral function from the target thread on the interval from the time of actual completion of the target operation until the estimated time point. Quantitatively absolute potential effect is determined on the unit time interval from the actual completion of the target operation by half of the product of the value added (cost) by the square of the estimated time interval.

In the paper, the problem of deriving a single interdisciplinary efficiency indicator (efficiency formula) was solved. This indicator is the ratio of the absolute potential effect to the resource intensity of the target operation. As a management criterion, this indicator is intended to be used as an optimization criterion in managed systems at all hierarchical levels.

Using the efficiency criterion as the single optimization criterion of all systems at all hierarchical levels opens up prospects for the full automation of management processes, taking into account economically sound decision–making.